\documentclass[12pt]{amsart}
\usepackage{lipsum}
\usepackage{amssymb}
\usepackage{enumerate}
\usepackage{amsthm}
\usepackage{tikz-cd}
\usepackage{amsmath}
\usepackage[mathscr]{euscript}
\usepackage[utf8]{inputenc}
\usepackage[english]{babel}
\usepackage{hyperref}
\hypersetup{
    colorlinks=true,
    linkcolor=blue,
    citecolor=red,
    filecolor=red,
    urlcolor=violet,
}

\makeatletter
\@namedef{subjclassname@2020}{%
  \textup{2020} Mathematics Subject Classification}
\makeatother

\newtheorem{thm}{Theorem}[section]
\newtheorem{lem}[thm]{Lemma}
\newtheorem{prop}[thm]{Proposition}
\newtheorem{defi}[thm]{Definition}

\newtheorem{xrem}[thm]{Remark}

\DeclareMathOperator{\Nef}{{Nef}}
\DeclareMathOperator{\mult}{{mult}}

\DeclareMathOperator{\Amp}{{Amp}}

\DeclareMathOperator{\Eff}{{Eff}}

\DeclareMathOperator{\Div}{{Div}}

\DeclareMathOperator{\vol}{{vol}}

\begin{document}
\baselineskip=15pt
\subjclass[2020]{Primary 14C20; Secondary 14C25, 14E05 }
\keywords{$b$-divisors, birational geometry, ample, nef, big, pseudo-effective}

\author{Snehajit Misra}
\author{Nabanita Ray}

\address{Department of Mathematics and Computing, Indian Institute of Technology (Indian School of Mines) IIT (ISM) Dhanbad - 826004, Jharkhand, India.}
\email[Snehajit Misra]{misra08@gmail.com}

\address{Indraprastha Institute of Information Technology, Delhi, Okhla Industrial Estate, Phase
III, New Delhi, Delhi 110020, India.}
\email[Nabanita Ray]{nabanita@iiitd.ac.in}

\begin{abstract}
 In this article, we define the notion of ample Cartier $b$-divisor classes by using the notion of Seshadri constants for Cartier $b$-divisor classes. In particular, we have shown that the set of all ample Cartier $b$-divisor classes forms a convex cone inside the nef cone of Cartier $b$-divisor classes. Furthermore, we have studied various properties of these Cartier ample $b$-divisor classes. We have also given an equivalent characterization of big Cartier $b$-divisor classes in terms of volume function of the pseudo-effective Cartier $b$-divisor classes.  More specifically, we prove that the set of all big Cartier $b$-divisor classes form a convex cone.  Finally we have investigated how the nef Cartier $b$-divisor classes behave under the pullback.
\end{abstract}

\title{Positive cones of $b$-divisor classes}
\maketitle
\section{Introduction}
  The Riemann-Zariski space is introduced by Zariski \cite{Z44} as a tool for de-singularizing algebraic varieties (see \cite{ZS75} and \cite{V97} for more details on the definition  of  Riemann-Zariski space and history of its origin). The Riemann-Zariski space is introduced as a projective limit of all bi-rational morphisms over a projective space $X$. Subsequently, the notion of $b$-divisors on Riemann-Zariski spaces was introduced by Shokurov in his paper \cite{S1}. It provides the right setting to construct the Zariski decomposition of $b$-divisor classes (cf. \cite{BFJ09},\cite{KM}) and that of higher dimensional cycle classes (cf. \cite{FL}). Intersection theory of $b$-divisors are studied in \cite{DF1}. Many of its applications are found in \cite{S2},\cite{A05},\cite{C07},\cite{F12},\cite{BFF12},\cite{HL},\cite{DF1},\cite{DF2}, \cite{BJ22} and \cite{T24}. For example, generalized $b$-divisors are introduced in \cite{T24} to prove relative Yau-Tian Donaldson conjecture over a smooth Fano variety which suggest a deep link between differential and algebraic geometry. 

  One can define in a natural way a Weil and a Cartier divisor classes for $b$-divisors as projective and inductive limits of Weil and Cartier divisors on the bi-rational spaces modelled over a smooth projective variety $X$ (see Subsection 3.2). This construction is due to Shoukurov \cite{S2}.  In this article, we follow the terminology of \cite{DF1} and that of \cite{DF2}.
 Let $X$ be a smooth projective variety of dimension $n$ over an algebraically closed field $k$ of characteristics 0. The Riemann-Zariski space of $X$ is denoted by $\mathcal{X}$. A Weil
 $b$-divisor class $\alpha$ is a family of classes $\alpha_{X'}$ in real Neron-Severi groups $N^1(X')$ that are compatible under pushforward maps, where $X'$ runs over all smooth bi-rational models lying over $X$. The set of all Weil $b$-divisor classes forms an infinite-dimensional vector space, denoted by $w\text{-}N^1(\mathcal{X})$. It can be identified with the projective limit of the real Neron-Severi groups of all smooth models over $X$ (cf. Section 4.1). Similarly, one can describe the Cartier $b$- divisor classes, denoted by $c\text{-}N^1(\mathcal{X})$ as inductive limits (cf. Section 4.1).  Moreover, $c\text{-}N^1(\mathcal{X})$  is dense in $w\text{-}N^1(\mathcal{X})$. In \cite{DF1}, the authors prove that for any nef Weil $b$-divisor class $\alpha$, there always exists a net of nef Cartier $b$-divisor classes which converge to $\alpha$.

 In this article, we define the notion of ample Cartier $b$-divisor classes by introducing the notion of Seshadri constants for Cartier $b$-divisor classes. In particular, we have shown that the set of all ample Cartier $b$-divisor classes forms a convex cone inside the nef cone of Cartier $b$-divisor classes. Furthermore, we have studied various properties of these Cartier ample $b$-divisor classes (cf. Section 5). 

 Big Cartier $b$-divisor classes and its volume are described in \cite{T24}. In section 6, we have given an equivalent characterization of big Cartier $b$-divisor classes in terms of volume function of the pseudo-effective Cartier $b$-divisor classes. More specifically, we prove that the set of all big Cartier $b$-divisor classes form a convex cone inside $c\text{-}N^1(\mathcal{X}),$ denoted by $c\text{-}\text{Big}^1(\mathcal{X}).$ We have also exhibited various properties of  $c\text{-}\text{Big}^1(\mathcal{X})$ in Section 6.

 In Section 7, we consider the fact that morphism between projective varieties induces a  map between  the corresponding set of Cartier $b$-divisor classes \cite{BFF12}. We investigate how the nef or ample Cartier $b$-divisor classes behave under the pullback with respect to these maps.


 \section{Preliminaries}
 Throughout this article, all the algebraic varieties are assumed to be irreducible and defined over an algebraically closed  field of characteristic zero. In this section we recall the definition of various positive cones and their main properties.
 \subsection{Nef cone and Pseudo-effective cone}
 Let $X$ be a smooth projective variety of dimension $n$ and $\mathcal{Z}_k(X)$ (respectively $\mathcal{Z}^k(X)$) in $ \Div(X)$ denotes the free abelian group generated by $k$-dimensional (respectively $k$-codimensional) subvarieties on $X$. Given two cycles $\alpha,\beta$ of complimentary dimensions in $X$, we denote their intersection product denoted by $\alpha\cdot \beta \in \mathbb{Z}.$
 Two cycles $Z_1$,$Z_2 \in \mathcal{Z}_k(X)$ are said to be numerically equivalent, denoted by $Z_1\equiv Z_2$ if $Z_1\cdot \gamma =  Z_2\cdot\gamma $
  for all $\gamma \in \mathcal{Z}^k(X)$. The \it numerical groups \rm  $ N_k(X)_{\mathbb{R}}$ are defined as the quotient of
  $\mathcal{Z}_k(X) \otimes \mathbb{R}$ modulo numerical equivalence, and $N^k(X)_{\mathbb{R}}:= N_{n-k}(X)_{\mathbb{R}}$. The intersection product induces a perfect pairing
 \begin{align*}
N^k(X)_{\mathbb{R}} \times N_k(X)_{\mathbb{R}} \longrightarrow \mathbb{R}
\end{align*}
which implies $N^k(X)_{\mathbb{R}} \cong (N_k(X)_{\mathbb{R}})^\vee$. The convex cone generated by the set of all effective $k$-cycles in $N_k(X)_\mathbb{R}$ is denoted by $\Eff_k(X)$ and its closure $\overline{\Eff}_k(X)$ is called the \it pseudo-effective cone \rm  of $k$-cycles in $X$. We will use the notation $\alpha \geq 0$ to say that the numerical class $\alpha $ is pseudo-effective. For two numerical classes $\alpha$ and $\beta$, we denote  $\alpha \geq \beta$  if and only if $\alpha-\beta$ is pseudo-effective. The \it nef cone \rm of co-dimension $k$-cycles in $N^k(X)_{\mathbb{R}}$, denoted by $\Nef^k(X)$ are defined as follows:
\begin{align*}
 \Nef^k(X) := \Bigl\{ \alpha \in N^1(X) \mid \alpha \cdot \beta \geq 0 \hspace{2mm} \forall \beta \in \overline{\Eff}_k(X)\Bigr\}.
\end{align*}
A line bundle $L$ on $X$ is said to be \it ample \rm if some of its positive power correspond to an embedding of $X$ into some projective space, i.e. $$\phi_{L^{\otimes m}}: X\rightarrow \mathbb{P}(H^0(X, L^{\otimes m}))$$
is an embedding for some $m\in\mathbb{N}$. Thus, the set of all ample classes in $N^1(X)_{\mathbb{R}}$ forms a cone inside $\Nef^1(X)$, known as the ample cone of $X$, and it is denoted by $\Amp^1(X)$. Infact, we have the following results :
\begin{center}
    $\overline{\Amp^1(X)} = \Nef^1(X)$ and int$(\Nef^1(X)) = \Amp^1(X)$.
\end{center}

For a line bundle $L$ on $X$, we define the semigroup $$N(X,L) := \Bigl\{m \in \mathbb{Z}_{\geq 0} \mid H^0(X,L^{\otimes m})\neq 0\Bigr\}.$$
A line bundle $L$ on $X$ is called big if and only if $$\phi_m : X\dashrightarrow \mathbb{P}\bigl(H^0(X,L^{\otimes m})\bigr)$$ is bi-rational onto its image for some positive $m>0$ in $N(X,L)$. A divisor $D$ is called big if the corresponding line bundle $\mathcal{O}_X(D)$ is big. Moreover, $D$ is big if and only if there is a constant $C>0$ such that $$h^0\bigl(X,\mathcal{O}_X(mD)\bigr) \geq C\cdot m^n$$
for all sufficiently large $m\in N(X,L)$ (\cite{L1}, Lemma 2.2.3). The convex cone generated by big divisors in $N^1(X)$ is called \it big cone \rm of $X$, denoted by Big$^1(X)$.
In particular, the interior int$\bigl(\overline{\Eff}^1(X)\bigr)$ is the Big cone Big$^1(X)$ of big divisors in $N^1(X)_{\mathbb{R}}$, and $\overline{\text{Big}^1(X)} = \overline{\Eff}^1(X).$ 
We refer the
reader to \cite{L1},\cite{L2} for more details about these cones.
\vspace{2mm}

\subsection{Projection formula}
Let $Y$ be another smooth projective variety. For any regular morphism $\pi : X \longrightarrow Y$,  one defines the pushforward by $\pi$ of the cycle $[V]$ associated to any irreducible subvariety $V$ of $X$ by the following formula :
\begin{align*}
 & f_*[V] :=
  \begin{cases}    d[f(V)]\hspace{6mm} ; \text{if} \hspace{2mm}  \dim\bigl(f(V)) = \dim\bigl(V\bigr)\vspace{4mm} \\
  0\hspace{18mm}; \text{if} \hspace{2mm} \dim\bigl(f(V)\bigr) < \dim\bigl(V\bigr).
  \end{cases}
\end{align*}
where $d$ is the degree of the field extension $[K(V):K(f(V))]$. This map descends to a linear map $f_* : N_k(X)_{\mathbb{R}}\longrightarrow N_k(Y)_{\mathbb{R}}$.

Applying duality, we have the dual map $f^*:N^k(Y)\longrightarrow N^k(X)$ satifying the following projection formula : for any $\alpha \in N^k(Y)$ and $\beta\in N^l(X)$, we have

$$f_*(f^*\alpha\cdot \beta) = \alpha\cdot f_*\beta.$$
The negativity lemma states that for any proper birational morphism $\pi : X\longrightarrow Y$, and for any nef class $\alpha\in \Nef^1(X)$, we have $\pi^*\pi_*\alpha \geq \alpha$.

\subsection{Seshadri constants for nef line bundles}
The Seshadri constant of a nef line bundle $L \in \Nef^1(X)$ at $x \in X$ is defined as
\begin{align*}
 \varepsilon(X,L,x) := \inf_{x \in C} \hspace{1mm}   \Bigl\{\frac{L\cdot C}{\mult_xC}\Bigr\}
\end{align*}
where the infimum is taken over all irreducible curves in $X$ passing through $x$ having the multiplicity $\mult_xC$ at $x$. One can easily check that it is enough to take the infimum over irreducible and reduced curves $C$. For a nef line bundle $L$, the global Seshadri constant is
$$\varepsilon(X,L) := \inf\limits_{x\in X}\varepsilon(X,L,x).$$ 
One useful charecterisation of an ample line bundle using Seshadri constant is that a line bundle is ample if and only if its global Seshadri constant $\varepsilon(X,L) >0$. To know more about Seshadri constant on line bundle see \cite{PSC}.
    


The Seshadri constant $\varepsilon(X,L,x)$ of a nef divisor class $L$ at a point $x\in X$ is an interesting invariant of $L$ that measures local positivity around $x$ in several ways: some
numerical, some cohomological via asymptotic jet separation, and even via differential geometry or from arithmetic height theory. If $L\in \Amp^1(X)$ is an ample line bundle, then $\varepsilon(X,L,x) \leq \sqrt[n]{L^n}$ for all $x \in X$, where $n$ is the dimension of $X$ and $L^n$ is the  $n$ fold self-intersection of $L$. Hence, $\varepsilon(X,L,x) \in (0,\sqrt[n]{L^n}]$.

\subsection{Volume function for big line bundle }

 The volume of a line bundle  $L$ on $X$, denoted by $\vol_X(L)$, is the non-negative number defined as follows:
$$\vol_X(L):=\limsup\limits_{m \to \infty}\Bigl\{\frac{\dim_kH^0(X,L^{\otimes m})}{m^n/n!}\Bigr\}.$$
The limsup is actually a limit (see \cite{F},\cite{T}). The volume function is  $n$-homogeneous and invariant under numerical equivalence.  Hence we have a well-defined mapping from the N\'{e}ron-Severi group $N^1(X)_{\mathbb{Q}}$ to $\mathbb{Q}_{\geq 0 }$ which we again denote by $\vol_X$. In fact, the function $$\vol_X : N^1(X)_{\mathbb{Q}} \mapsto \mathbb{Q}_{\geq 0 }$$
$$\xi \mapsto \vol_X(\xi)$$ on $N^1(X)_{\mathbb{Q}}$ extends uniquely to a continuous function $\vol_X: N^1(X)_{\mathbb{R}} \longrightarrow \mathbb{R}_{\geq 0}$ on the real N\'{e}ron-Severi group $N^1(X)_{\mathbb{R}}$. There is another characterization of big line bundles in terms of volume function.  A divisor class $L$ in $N^1(X)_{\mathbb{R}}$ is big if and only if $\vol_X(L) >0$.
\vspace{1mm}

\section{Numerical $b$-cycle classes} 

\subsection{Riemann-Zariski space}
Let $X$ be a smooth projective variety of dimension $d$ over an algebraically closed field of characteristic 0. A smooth model over $X$ is a projective bi-rational morphism $\pi' : X'\longrightarrow X$ from a smooth projective variety $X'$.

We denote the category of all smooth models over $X$ by $\mathcal{M}_X$. By abuse of notation, we will consider $\mathcal{M}_X$ as a set consisting of all smooth varieties modeled over $X$. We define a partially ordered relation $\geq $ on $\mathcal{M}_X$ as follows : for any $X',X^{''}\in \mathcal{M}_X$, we say $X'\geq X^{''}$ if and only if the map $\pi := (\pi^{''})^{-1}\circ \pi' : X'\longrightarrow X^{''}$ is regular. Given any two elements $X_1,X_2 \in \mathcal{M}_X$, one can find $X_3\in \mathcal{M}_X$ such that $X_3\geq X_1$ and $X_3\geq X_2$. It is convenient to introduce the projective limit over the inductive  system $\mathcal{M}_X$ of all smooth models over $X$, each endowed with the Zariski topology. Let $$\mathcal{X}:= \varprojlim X'$$ be the projective limits over $\mathcal{M}_X$. One obtains in this way a quasi-compact topological space which is the Riemann-Zariski space attached to $X$. Note that when the base field $K$ is countable then it is separable and metrizable, but these properties do not hold when $K$ is
uncountable.
\subsection{Weil and Cartier $b$-cycle classes}\label{subsec4.1} A Weil $b$-cycle class $\alpha$ of co-dimension $k$ is a map which assign any smooth model $X' \in \mathcal{M}_X$ a numerical class $\alpha_{X'}\in N^k(X')$ such that for any pair of smooth models $X'\geq X''$ with $\pi := (\pi^{''})^{-1}\circ \pi' : X'\longrightarrow X^{''}$ regular, we have $\pi_*(\alpha_{X'})=\alpha_{X''}$. The numerical class $\alpha_{X'}$ is called the incarnation of $\alpha$ in $X'$. When the dimension of $X$ is at least 2, the space of Weil $b$-cycle classes of co-dimension $k$ is an infinite dimensional real vector space, which we denote by $w\text{-}N^k(\mathcal{X})$. The space $w\text{-}N^k(\mathcal{X})$ can be seen as the projective limit of the spaces $\bigl\{N^k(X')_{\mathbb{R}} \mid X'\in \mathcal{M}_X\bigr\}$, where morphisms are given by pushforwards on smooth models, i.e.
\begin{center}
$w\text{-}N^k(\mathcal{X}) := \Bigl\{(\alpha_{X'})_{X'\in \mathcal{M}_X} \in \prod\limits_{X'\in \mathcal{M}_X}N^k(X')_{\mathbb{R}}\mid \pi_*(\alpha_{X'}) = \alpha_{X''}\hspace{1mm} \text{for}\hspace{1mm}
X'\geq X^{''}\Bigr\}.$
\end{center}
It is thus endowed with a natural product (or weak) topology. Thus the space $w\text{-}N^k(\mathcal{X})$ is a locally convex topological vector space when endowed with its
natural product topology such that every bounded set is relatively compact.

A Cartier $b$-cycle class of co-dimension $k$ is a Weil $b$-cycle class $\alpha$ for which one can find a model $X''\in \mathcal{M}_X$ and a  divisor class $\beta \in N^k(X^{''})$ such that $\alpha_{X'} = \pi^*\beta$ for any smooth model $X'\geq X^{''}$, where $\pi := (\pi^{''})^{-1}\circ \pi' : X'\longrightarrow X^{''}$ is regular. When it is the case, we say $\alpha$ is determined in $X''$.
A determination of a Cartier $b$-divisor class $\alpha$ is a pair consisting of a smooth model $X''$ and a class $\beta\in N^k(X'')$ such that $\alpha$ is determined in $X''$ and $\alpha_{X''}=\beta$.
The set of Cartier $b$-cycle classes of co-dimension $k$ will be denoted by $c\text{-}N^k(\mathcal{X})$ as a subspace of  $w\text{-}N^k(\mathcal{X})$.  Note that  if $X''\geq X' \geq X_0,$ then $\alpha_{X''}=\pi^*\alpha_{X'},$ where $\pi:X''\longrightarrow X'$ is the regular map

One can see that the set of Cartier $b$-cycle classes as the inductive limits of $\{N^k(X')_{\mathbb{R}}\mid X'\in \mathcal{M}_X\}$, and thus $c\text{-}N^k(\mathcal{X})$ is a locally convex topological vector space with the inductive limit topology. More precisely, there is an equivalence relation on the set $$\bigcup\limits_{X'\in \mathcal{M}_X}N^k(X')$$ defined as follows :
 For $\alpha\in N^k(X')$ and $\beta\in N^k(X'')$ for some $X',X''\in \mathcal{M}_X$, we say $\alpha\sim  \beta$ if and only if there is another $\tilde{X}$ such that $\tilde{X}\geq X'$, $\tilde{X}\geq X''$ satifying the following diagram
 \begin{center}
\begin{tikzcd}[column sep=tiny]
& \tilde{X} \ar[dr, "\tilde{\pi}_2"] \ar[dl, "\tilde{\pi}_1"']
&
&[1.5em] \\
X_1   \ar[dr, "\pi_1"']
&
& X_2  \ar[dl, "\pi_2"] \\
& X
&
&
\end{tikzcd}
\end{center}

 and $\tilde{\pi}_1^*(\alpha)=\tilde{\pi}_2^*(\beta)$. Then $$c\text{-}N^k(\mathcal{X}) := \varinjlim N^k(X') =  \bigcup\limits_{X'\in \mathcal{M}_X} N^k(X')/\sim .$$

 We denote by $[\alpha]$ the Cartier $b$-class associated to a class $\alpha\in N^k(X^{'})_{\mathbb{R}}$.
For each smooth model $X'\in \mathcal{M}_X$, the map $$\alpha \mapsto [\alpha]$$ induces an injective linear map $$\phi_{X'}: N^k(X')_{\mathbb{R}}\longrightarrow c\text{-}N^k(\mathcal{X}).$$  In this topology, we have a set $U \subseteq c\text{-}N^k(\mathcal{X})$ is open if and only if $\phi_{X'}^{-1}(U)\subseteq N^k(X')_{\mathbb{R}}$ is open for each $\phi_{X'} : N^k(X')_{\mathbb{R}}\longrightarrow c\text{-}N^k(\mathcal{X})$.
The space $c\text{-}N^k(\mathcal{X})$ is dense in $w\text{-}N^k(\mathcal{X})$.

\begin{xrem}
   \rm Let $K_{\mathcal{X}}$ be a Weil $b$-divisor class on $X$ i.e. $K_{\mathcal{X}}\in w\text{-}N^1(\mathcal{X})$ such that $K_{{\mathcal{X}}_{X'}}=K_{X'}$, where $X'\in\mathcal{M}_X$ and $K_{X'}$ is the canonical divisor on $X'$. Then $K_{\mathcal{X}}$ is called the canonical $b$-divisor class\it.\rm
\end{xrem}

\begin{xrem}
\rm A Cartier divisor $D$ on a given model $X'\in\mathcal{M}_X$ induces a Cartier $b$-divisor $\alpha_D$. It is defined by pulling back $D$ to all models that dominate $X'$ and then pushing forward on all models. By definition, all Cartier $b$-divisor classes are actually obtained in this way.
\end{xrem}
\subsection{Intersection Product and a perfect pairing}
Let $\alpha\in c \text{-} N^k(\mathcal{X})$ which determined in a model $X_0\in \mathcal{M}_X$  and $\beta\in w \text{-} N^{l}(\mathcal{X})$. For any $X'\geq X_0$,  we set $$(\alpha\cdot\beta)_{X'}:=\alpha_{X'}\cdot\beta_{X'}= \pi^*\alpha_{X_0}\cdot \beta_{X'},$$ where $\pi : X'\longrightarrow X_0$ is the regular map. Recall that  if $X''\geq X' \geq X_0,$ then $\alpha_{X''}=\pi^*\alpha_{X'},$ where $\pi:X''\longrightarrow X'$ is the regular map. Now by projection formula $$\pi_*(\alpha\cdot\beta)_{X''}=\pi_*(\pi^*\alpha_{X'}\cdot\beta_{X''})=\alpha_{X'}\cdot \pi_*\beta_{X''}=\alpha_{X'}\cdot \beta_{X'}=(\alpha.\beta)_{X'}.$$ Therefore one can define the class $\alpha\cdot\beta\in w\text{-}N^{k+l}(\mathcal{X})$ as the unique Weil $b$-class whose incarnation in any smooth model $X'$ dominating $X_0$ is equal to $(\alpha\cdot \beta)_{X'}$. Note that for any smooth model $X'$ such that $X_0\geq X'$ so that the canonical birational map $\pi : X_0\longrightarrow X'$ is regular, the incarnation of $(\alpha\cdot \beta)_{X'} = \pi_*(\alpha\cdot \beta)_{X_0}$, which
is different from $(\alpha_{X'}\cdot \beta_{X'})$ in general. 

\begin{lem}
  The   intersection of two Cartier $b$-classes is again a Cartier $b$-class.
  \begin{proof}
      Let $\alpha$ and $\beta$ be two Cartier $b$-classes determined on $X'$ and $X''$ respectively. There exists $\tilde{X}\in\mathcal{M}_X$ such that $\tilde{X}\geq X', X''$. Therefore, $\alpha$ and $\beta$ both are determined on $\tilde{X}$. Enough to show $\tilde{X}$ is a determination of $\alpha\cdot\beta$. Let $X_0\geq \tilde{X}$ and $\pi: X_0\rightarrow\tilde{X}$ is the corresponding regular birational map. Hence $\pi^*(\alpha_{\tilde{X}}\cdot \beta_{\tilde{X}})=\pi^*\alpha_{\tilde{X}}\cdot\pi^*\beta_{\tilde{X}}=\alpha_{X_0}\cdot\beta_{X_0}=(\alpha.\beta)_{X_0}$.
  \end{proof}
\end{lem}

Let $w\text{-}N^*(\mathcal{X}) = \bigoplus\limits_{k}N^k(\mathcal{X})$ and $c\text{-}N^*(\mathcal{X}) = \bigoplus\limits_{k}N^k(\mathcal{X})$. Then the intersection product defined above induces a bilinear pairing
$$c\text{-}N^*(\mathcal{X}) \times w\text{-}N^{*}(\mathcal{X}) \longrightarrow w\text{-}N^{*}(\mathcal{X}),$$
respecting the natural grading, and continuous for the weak topology in the
second variable.

We observe that $c\text{-}N^n(\mathcal{X}) = w\text{-}N^n(\mathcal{X})$ is canonically isomorphic to $\mathbb{R}$ so that
Poincar\'{e} duality on each smooth model implies that the induced pairing
$$c\text{-}N^k(\mathcal{X}) \times w\text{-}N^{d-k}(\mathcal{X}) \longrightarrow \mathbb{R}$$ is a perfect pairing.
\subsection{Nef and pseudo-effective $b$-classes}
The convex cone in $c\text{-}N^1(\mathcal{X})$ generated by Cartier numerical $b$-divisor class $[\alpha]$ with $\alpha\in \Nef^1(X')$ for some smooth model $X'\in \mathcal{M}_X$ is denoted by $c\text{-}\Nef^1(\mathcal{X})$.  A class $\alpha \in c\text{-}N^1(\mathcal{X})$ is called nef if $\alpha\in c\text{-}\Nef^1(\mathcal{X}).$
Since nef classes are stable by pull-back, a class $\alpha\in c\text{-}N^1(\mathcal{X})$ is nef if and only if it is nef in one (or any) of its determination.   For any $\alpha\in c\text{-}N^1(\mathcal{X}),$  if $\alpha_{X_0} \in \Nef^1(X_0)$ in one of its determination $X_0$, then it is nef in any of its determination. Also,
s divisor class $\alpha \in c\text{-}N^1(\mathcal{X})$ is nef if and only if it is nef in one (or any) of its determination.    

\begin{defi}
 The weak closure of $c\text{-}\Nef^1(\mathcal{X})$ in $w\text{-}N^1(\mathcal{X})$ is called the nef cone of $b$-divisor classes, and it is denoted by $\Nef^1(\mathcal{X})$. Hence $$\Nef^1(\mathcal{X}) = w\text{-}N^1(\mathcal{X})\cap c\text{-}\Nef^1(\mathcal{X}).$$
\end{defi}
A class $\alpha\in w\text{-}N^1(\mathcal{X})$ is nef if and only if there exists a net of $b$-classes $\alpha_i\in c\text{-}\Nef^1(\mathcal{X})$ such that $(\alpha_i)_{X'}\longrightarrow \alpha_{X'}$ for all smooth model $X'\in \mathcal{M}_{X}$.

\begin{defi}
 A class $\alpha\in w\text{-}N^k(\mathcal{X})$ is said to be pseudoeffective if $\alpha_{X'}\geq 0$ for all smooth models $X'\in \mathcal{M}_X$.
\end{defi}
 Note that pseudo-effectivity is only preserved by pushforwards, and not by pullback, so that it may happen that a Cartier $b$-class determined by a pseudoeffective class in a smooth model $X$ is not pseudoeffective in $w\text{-}N^k(\mathcal{X})$. But in codimension 1, a Cartier $b$-class $\alpha \in c\text{-}N^1(\mathcal{X})$ is pseudoeffective if and only if one (or any) of its determination on a smooth model is pseudoeffective.


\section{Ample $b$-divisor classes}\label{sec5}
In this section, we introduce the notion of ample $b$-divisor classes and the notion of ample cone of $b$-divisor classes. We also exhibit some topological properties of these ample cones in $ c\text{-}N^1(\mathcal{X})$.

\begin{defi}
 Let $\alpha $ be a  Cartier $b$-divisor nef class in $c\text{-}\Nef^1(\mathcal{X})$. Then the Seshadri constant of the $b$-divisor class $\alpha$, denoted by $\varepsilon(\mathcal{X},\alpha)$ is defined as follows :
 \begin{align}\label{seq1}
 \varepsilon(\mathcal{X},\alpha) : = \sup\limits_{X'\in \mathcal{M}_X}\Bigl\{\varepsilon(X',\alpha_{X'})\Bigr\},
 \end{align}
 where the supremum is taken over all determination $\alpha_{X'}$ and for all $X'\in \mathcal{M}_X$.
 \end{defi}
Let $X_1, X_2$ be two determination of $\alpha \in c\text{-}N^1(\mathcal{X})$, then they may not be comparable. But there always exists $\tilde{X}\geq X_1$ and $\tilde{X}\geq X_2$ such that
\begin{center}
\begin{tikzcd}[column sep=tiny]
& \tilde{X} \ar[dr, "f_2"] \ar[dl, "f_1"']
&
&[1.5em] \\
X_1   \ar[dr, "\pi_1"']
&
& X_2  \ar[dl, "\pi_2"] \\
& X
&
&
\end{tikzcd}
\end{center}
with
\begin{align*}
&\varepsilon(X_1,\alpha_{X_1})\leq \sqrt[n]{\alpha_{X_1}^n}\\
&\varepsilon(X_2,\alpha_{X_2})\leq \sqrt[n]{\alpha_{X_2}^n},
\end{align*}
where $n$ is the dimension of $X$. Observe that
\begin{align*}
& \alpha_{X_1}^n = \deg(f_1)(f_1^*\alpha_{X_1})^n = \alpha_{\tilde{X}}^n\\
& \alpha_{X_2}^n = \deg(f_2)(f_2^*\alpha_{X_2})^n = \alpha_{\tilde{X}}^n.
\end{align*}
This shows that 
\begin{center}
$\Bigl\{\varepsilon(X',\alpha_{X'})\mid X'\in \mathcal{M}_X, X'$ is a determination of $\alpha \Bigr\}$
\end{center}
is bounded above. Thus the supremum exists in (\ref{seq1}).
\begin{lem}
    For a Cartier $b$-divisor class  $\alpha\in c\text{-}N^1(\mathcal{X})$, we have $\alpha\in c\text{-}\Nef^1(\mathcal{X})$ if and only if $\varepsilon(\mathcal{X},\alpha) \geq 0.$
    \begin{proof}
 As Seshadri constant of a nef divisor on a projective variety is non-negative, $\varepsilon(\mathcal{X},\alpha) \geq 0$ follows from the definition.  
    \end{proof}
\end{lem}
\begin{defi}
 A Cartier $b$-divisor class $\alpha\in c\text{-}\Nef^1(\mathcal{X})$ is said to be ample $b$-divisor class if $\varepsilon(\mathcal{X},\alpha)>0$. The set of all ample Cartier $b$-divisor classes will be denoted by $c\text{-}\Amp^1(\mathcal{X})$.
\end{defi}
Note that by the definition of Seshadri constant for $b$-divisor classes $$c\text{-}\Amp^1(\mathcal{X}) \subset c\text{-}\Nef^1(\mathcal{X}).$$
\begin{defi}
 Let $\alpha,\beta \in c\text{-}N^1(\mathcal{X})$ be two Cartier $b\text{-}$divisor classes and $\lambda\in \mathbb{R}_{>0}$ be any positive real number. Then we define the numerical $b$-divisor classes $\alpha+\beta$ and $\lambda\alpha$ as follows : for any smooth model $X'\in \mathcal{M}_X$, we define $$(\alpha+\beta)_{X'} := \alpha_{X'}+\beta_{X'}.$$ and $$ (\lambda\alpha)_{X'} := \lambda\alpha_{X'}.$$
 \end{defi}
 Let $\alpha$ is determined by $X_0$. Then for any $X'\geq X_0$ and the regular map $\pi : X'\longrightarrow X_0$, we have $$ \pi^*((\alpha+\beta)_{X_0}) = \pi^*(\alpha_{X_0}+\beta_{X_0}) = \pi^*(\alpha_{X_0})+\pi^*(\beta_{X_0}) = \alpha_{X'}+\pi^*\beta_{X_0}.$$
\begin{thm} \label{thm5.5}
 The set $c\text{-}\Amp^1(\mathcal{X})$ forms a convex cone inside $c\text{-}\Nef^1(\mathcal{X})$.
 \begin{proof}
Suppose $\alpha,\beta \in c\text{-}\Amp^1(\mathcal{X})$ such that $\varepsilon(\mathcal{X},\alpha)>0$ and $\varepsilon(\mathcal{X}, \beta)>0$. Then there exist smooth models $X',X''\in \mathcal{M}_X$ such that $\varepsilon(X',\alpha_{X'})>0$ and $\varepsilon(X'',\beta_{X''})>0$. Therefore $\alpha_{X'}$ and $\beta_{X''}$ are ample divisors on $X'$ and $X''$ respectively. Note that $Y=  X' \times_{X} X'' \in \mathcal{M}_X$. Now consider the following diagram
\begin{center}
\begin{tikzcd}
Y=  X' \times_{X} X'' 
\arrow[drr, bend left, "\phi''"]
\arrow[ddr, bend right, "\phi'"]
\arrow[dr,  "\phi"] & & \\
& X' \times_{Spec(k)} X'' \arrow[r, "\pi''"] \arrow[d, "\pi'"]
& X'' \arrow[d] \\
& X' \arrow[r]
& X \arrow[r]
& Spec(k)
\end{tikzcd}
\end{center}
The map $\phi$ in the above diagram is a closed embedding.
Now, observe that $\pi'^{*}(\alpha_{X'})+\pi''^{*}(\beta_{X''})$ is an ample divisor on $X' \times_{Spec(k)} X''$. Therefore $\phi'^{*}(\alpha_{X'})+\phi''^{*}(\beta_{X''})=\alpha_Y+\beta_Y$ is also an ample divisor on $Y=  X' \times_{X} X'' $, as the restriction of ample divisor is also an ample divisor. Now $ \varepsilon(Y,\alpha_Y+\beta_Y) >0$. Hence
\begin{align*}
 \varepsilon(\mathcal{X},\alpha+\beta)  >\varepsilon(Y,\alpha_Y+\beta_Y) > 0.
\end{align*}
Similarly, for any positive $\lambda \in \mathbb{R}_{> 0}$, we have
\begin{align*}
 \varepsilon(\mathcal{X},\lambda\alpha) &
= \sup\limits_{X'\in \mathcal{M}_X} \Bigl\{\varepsilon(X',(\lambda\alpha)_{X'}\Bigr\}  \\
& = \sup\limits_{X'\in \mathcal{M}_X} \Bigl\{\varepsilon(X',\lambda\cdot \alpha_{X'})\Bigr\} \\
& = \lambda \cdot \varepsilon(\mathcal{X}, \alpha) > 0
\end{align*}
 This shows that the set $c\text{-}\Amp^1(\mathcal{X})$ forms a convex cone inside $c\text{-}N^1(\mathcal{X})$.
 \end{proof}

\end{thm}
Let us recall the injective linear map $\phi_X': N^k(X')_{\mathbb{R}}\longrightarrow c\text{-}N^k(\mathcal{X})$, for any smooth model $X'\in\mathcal{M}_X$ given by $\phi_{X'}(\alpha)=[\alpha]$. Note that this Cartier b-class $[\alpha]$ is determined by $X'$.
\begin{thm}
 The set $c\text{-}\Nef^1(\mathcal{X})$ is closed with respect to inductive topology in $c\text{-}N^1(\mathcal{X})$.
 \begin{proof}
It is enough to show that $U:= c\text{-}N^1(\mathcal{X})\setminus c\text{-}\Nef^1(\mathcal{X})$ is open subset of $c\text{-}N^1(\mathcal{X})$ with respect to inductive topology. We will show that $\phi_{X'}^{-1}(U)$ is open for every $X'\in \mathcal{M}_X$ for the natural injective map $\phi_{X'} : N^1(X')\longrightarrow c\text{-}N^1(\mathcal{X})$.

Note that
\begin{align*}
\alpha \in \phi_{X'}^{-1}(U) &  \iff [\alpha] \text{ determined by $X'$ and } \phi_{X'}(\alpha) = [ \alpha ] \notin c\text{-}\Nef^1(\mathcal{X})\\
& \iff [\phi_{X'}(\alpha)]_{X'} = \alpha \notin \Nef^1(X').
\end{align*}
This shows that $\phi_{X'}^{-1}(U) = N^1(X')\setminus \Nef^1(X')$ for every $X'\in \mathcal{M}_X$. Hence the result follows.
 \end{proof}

\end{thm}


\begin{thm}
 The closure of the convex cone $c\text{-}\Amp^1(\mathcal{X})$ is equal to  $c\text{-}\Nef^1(\mathcal{X})$.
 \begin{proof}
Let $\alpha \in c\text{-}\Nef^1(\mathcal{X})$. We will show that every open set $U$ containing $\alpha$ satisfies $$U\cap c\text{-}\Amp^1(\mathcal{X}) \neq \emptyset.$$
Let us assume that $\alpha$ is determined by $X'\in\mathcal{M}_X$. Hence $\alpha_{X'}\in \Nef^1(X')$.
Note that $ \phi_{X'}^{-1}(U)\subseteq N^1(X')$ is an open neighborhood of $\alpha_{X'}$. Hence we get
\begin{align*}
 \phi_{X'}^{-1}(U)\cap \Amp^1(X') \neq \emptyset & \implies U\cap \phi_{X'}(\Amp^1(X')) \neq \emptyset\\
 & \implies U\cap c\text{-}\Amp^1(\mathcal{X}) \neq \emptyset, \hspace{2mm}
 \text{as}\hspace{2mm} \phi_{X'}(\Amp^1(X'))\subseteq c\text{-}\Amp^1(\mathcal{X})
\end{align*}
 \end{proof}
\end{thm}
 \begin{defi}
 A Weil class $\alpha\in w\text{-}N^k(\mathcal{X})$ is said to be almost zero if its any incarnation is zero i.e. $\alpha_{X'}=0$ for some $\alpha_{X'}\in \mathcal{M}_X$. We denote the set of all most zero classes by $\mathcal{Z}_X$.
 \end{defi}
For any $\alpha\in c\text{-}N^1(\mathcal{X})$, we set

$$\alpha_{> 0}:=\big\{\beta\in w\text{-}N^{d-1}(\mathcal{X}) \mid \alpha\cdot\beta>0\big\}.$$

\begin{prop}
Let $\alpha$ be a Cartier class. If $\alpha\in c\text{-Amp}^1(\mathcal{X})$, then 
$$\overline{{\Eff}}^{d-1}(\mathcal{X})\setminus \mathcal{Z}_X\subset \alpha_{> 0}.$$
\end{prop}
\begin{proof}
 As $\alpha\in c\text{-Amp}^1(\mathcal{X})$, there exists a determination $\alpha_{X'}$ such that  $0<\varepsilon(\mathcal{X},\alpha )\leq\varepsilon(X',\alpha_{X'} )$. Therefore $\alpha_{X'}$ is an ample line bundle on $X'\in \mathcal{M}_X$. So for any $\beta\in\overline{{\Eff}}^{d-1}(\mathcal{X})\setminus \mathcal{Z}_X$ we have $\alpha\cdot\beta=\alpha_{X'}\cdot\beta_{X'}>0$. Hence, the result follows.
\end{proof}
The above result produces the necessary condition for $\alpha $ to be an ample Cartier $b\text{-}$ classes. Another result in this direction is the following.

\begin{thm}
Let $\alpha$ be a Cartier b-class. Then $\alpha\in c\text{-Amp}^1(\mathcal{X})$ if and only if there exists $X'\in\mathcal{M}_X$ such that 
$$\phi_{X'}\big(\overline{{\Eff}}^{d-1}(X')\setminus\{0\}\big)\subset \alpha_{>0},$$ where
$\phi_{X'}: N^{d-1}(X')\rightarrow c\text{-}N^{d-1}(\mathcal{X})$.
\end{thm}
\begin{proof}
This result is immediate  from the fact that a divisor $D$ on a projective variety $Y$ is ample if and only if $\overline{{\Eff}}^{d-1}(Y)\setminus\{0\}\subseteq D_{>0}$ (see \cite{L1} Theorem 1.4.29).
\end{proof}
Note that $\phi_{X'}\big(\overline{{\Eff}}^{d-1}(X')\setminus\{0\}\big)$ may not be contained in $\overline{{\Eff}}^{d-1}(\mathcal{X})$, because pull back of pseudo-effective classes may not be pseudo-effective.
\section{Big $b$-divisor classes}\label{sec6}
Recall that a line bundle $L$ on a smooth projective variety $X$ is big if and only if there exists a constant $C>0$ such that $$h^0(X,L^{\otimes m}) \geq Cm^n$$
for all sufficiently large $m\in N(X,l) := \Bigl\{ m\in \mathbb{N} \mid h^0(X,L^{\otimes m})\neq 0\Bigr\}.$
\begin{lem}
    Let $X',X''\in \mathcal{M}_X$ be two smooth projective varieties over $X$ and $\pi : X^{''}\longrightarrow X'$ is a proper bi-rational map. Further assume that $L$ is a big line bundle on $X'$. Then $\pi^*L$ is also big line bundle on $X''.$
    \begin{proof}
    Since $L$ is big, there exists a constant $C>0$ such that $$h^0(X',L^{\otimes m}) \geq Cm^n$$ for all sufficiently large $m\in N(X',L).$

    Note that by projection formula
    \begin{align*}
         h^0(X'',\pi^*L^{\otimes m}) 
        = h^0(X',\pi_*\mathcal{O}_{X''}\otimes L^{\otimes m}) 
        = h^0(X',L^{\otimes m}) 
    \end{align*}
    This implies $$ h^0(X'',\pi^*L^{\otimes m}) \geq C\cdot m^n$$ for sufficiently large $m\in N(X'',\pi^*L).$ In other words, $\pi^*L$ is big. 
    \end{proof}
\end{lem}
Note that a Cartier $b$-divisor $\alpha\in c\text{-} N^1(\mathcal{X})$ is pseudo-effective if and only if its incarnation in one (or all) of its determination on a smooth model is pseudo-effective.
So for a $b$-divisor class $\alpha\in c\text{-}\overline{{\Eff}}^1(\mathcal{X}),$ we define the volume function
$$c\text{-}\vol_{\mathcal{X}}(\alpha) := \inf\limits_{X'\in \mathcal{M}_X} \Bigl\{\vol(\alpha_{X'})\Bigr\},$$
where the infimum is taken over all determination of $\alpha$ of smooth models over $X$.
\begin{xrem}
Let $\alpha\in c\text{-}\overline{{\Eff}}^1(\mathcal{X})$ and $\alpha$ be determined on $X'$ and $X''$. We know that there exists $\tilde{X}\in\mathcal{M}_X$ such that

\begin{center}
\begin{tikzcd}[column sep=tiny]
& \tilde{X} \ar[dr, "f_2"] \ar[dl, "f_1"']
&
&[1.5em] \\
X'  \ar[dr, "\pi_1"']
&
& X''  \ar[dl, "\pi_2"] \\
& X
&
&
\end{tikzcd}
\end{center}
Note that $\alpha_{\tilde{X}}=f_1^*\alpha_{X'}=f_2^*\alpha_{X''}$. Therefore, by the projection formula $h^0(X',\alpha_{X'})=h^0(X'',\alpha_{X''})=h^0(\tilde{X},\alpha_{\tilde{X}}).$ So clearly $\vol(\alpha_{X'})=\vol(\alpha_{X''})=\vol(\alpha_{\tilde{X}})$. Hence the volume of any Cartier classes is same as the the volume of any (or all) of its determination i.e. $c\text{-}\vol_{\mathcal{X}}(\alpha)=\vol(\alpha_{X'})$, where $X'$ is a determination of $\alpha$.
\end{xrem}

\begin{defi}
 A $b$-divisor class $\alpha\in c\text{-}\overline{{\Eff}}^1(\mathcal{X})$ is said to be c-big if $c\text{-}\vol_{\mathcal{X}}(\alpha)>0$.
\end{defi}

\begin{thm} \label{thm5.5}
 The set $c\text{-}\text{Big}^1(\mathcal{X})$ forms a convex cone inside $c\text{-}N^1(\mathcal{X})$.
 \begin{proof}
Suppose $\alpha,\beta \in c\text{-}\text{Big}^1(\mathcal{X})$ such that $c\text{-}\vol_{\mathcal{X}}(\alpha)>0$ and $c\text{-}\vol_{\mathcal{X}}(\beta)>0$. Consider that $\alpha$ and $\beta$ are determined  on smooth models $X',X''\in \mathcal{M}_X$. Then $\vol(\alpha_{X'})>0$ and $\vol(\beta_{X''})>0$. Therefore $\alpha_{X'}$ and $\beta_{X''}$ are big divisors on $X'$ and $X''$ respectively.
There exists a smooth model $\tilde{X}\geq X_1$ and $\tilde{X}\geq X_2$ such that $\alpha, \beta $ and $\alpha+\beta$ all are determined on $\tilde{X}$. As pull back of big divisors are big,  $f_1^{*}(\alpha_{X'})+f_2^{*}(\beta_{X''})=\alpha_{\tilde{X} }+\beta_{\tilde{X} }$ is a big divisor on $\tilde{X}$. Hence $c\text{-}\vol_{\mathcal{X}}(\alpha+\beta)=\vol(\alpha_{\tilde{X} }+\beta_{\tilde{X} })>0$

Similarly, for any positive $\lambda \in \mathbb{R}_{> 0}$, we have $\lambda\alpha$ is also a Cartier class.

 \end{proof}

\end{thm}

\begin{thm}
 The convex cone $c\text{-}\text{\rm Big}^1(\mathcal{X})$ is open in $c\text{-}N^1(\mathcal{X})$.
 \begin{proof}
 We will show that for all $X'\in \mathcal{M}_X$ and for the injective map $$\phi_{X'} : N^1(X')\longrightarrow c\text{-}N^1(\mathcal{X}),$$ one has $$\phi_{X'}^{-1}\bigl(c\text{-}\text{Big}^1(\mathcal{X})\bigr) = \text{Big}(X').$$
 
Clearly, $$\text{Big}(X') \subseteq \phi_{X'}^{-1}\bigl(c\text{-}\text{Big}^1(\mathcal{X})\bigr).$$
Now, let $\alpha\in c\text{-}\text{Big}^1(\mathcal{X}).$ Then $c\text{-}\vol(\alpha) > 0.$ This implies
\begin{align*}
   c\text{-}\vol(\alpha) = \inf\limits_{X'\in \mathcal{M}_X} \vol(\alpha_{X'}) > 0
    & \implies \vol(\alpha_{X'})>0 \hspace{2mm} \text{for} \hspace{1mm} \text{all} \hspace{1mm} X'\in \mathcal{M}_{X} \text{ and $\alpha$ is determined on $X'$ },\\
    & \implies \alpha_{X'} \in \text{Big}(X')
\end{align*}
This shows that $$\phi_{X'}^{-1}\bigl(c\text{-}\text{Big}^1(\mathcal{X})\bigr) \subseteq \text{Big}(X').$$ This completes the proof.
 \end{proof}

\end{thm}
\begin{thm}
The interior of $c\text{-}\overline{{\Eff}}^1(\mathcal{X})$  is $c\text{-}\text{\rm Big}^1(\mathcal{X})$.
\begin{proof}
Clearly $c\text{-}\text{\rm Big}^1(\mathcal{X}) \subseteq \text{int}\Bigl( c\text{-}\overline{{\Eff}}^1(\mathcal{X})\Bigr).$

To prove the converse, let $U:= \text{int}\Bigl( c\text{-}\overline{{\Eff}}^1(\mathcal{X})\Bigr).$
For all $X'\in \mathcal{M}_X$, and for the injective map $$ \phi_{X'} : N^1(X')\longrightarrow c\text{-}N^1(\mathcal{X}), $$ one has $ \phi_{X'}^{-1}(U)$  is an open subset of $\phi_{X'}^{-1}\bigl(\overline{\Eff}^1(\mathcal{X})\Bigr) = \overline{\Eff}^1(X').$ 

Now $\text{Big}^1(X')$ is the interior of  $\overline{\Eff}^1(X')$. Thus $\text{Big}^1(X')$ is the largest open set contained in $\overline{\Eff}^1(X')$. Therefore we conclude that $$\phi_{X'}^{-1}(U) \subseteq \text{Big}(X').$$

Now let $\alpha \in U$ so that $\alpha_{X'}\in \phi_{X'}^{-1}(U).$ So $\vol(\alpha_{X'})>0$ for all $X'\in \mathcal{M}_X.$
    Thus $c\text{-}\vol(\alpha) > 0,$ or in other words $\alpha\in c\text{-}\text{Big}^1(\mathcal{X}).$ This implies $$U\subseteq c\text{-}\text{Big}^1(\mathcal{X}).$$
    This completes the proof.
\end{proof}
\end{thm}
Note that $\phi_{X'}^{-1}\bigl(\overline{\Eff}^1(\mathcal{X})\Bigr) = \overline{\Eff}^1(X')$  for all $X'\in \mathcal{M}_X$ where $ \phi_{X'} : N^1(X')\longrightarrow c\text{-}N^1(\mathcal{X}).$ So by the inductive limit topology $c\text{-}\overline{{\Eff}}^1(\mathcal{X})$ is a closed subset of $ c\text{-}N^1(\mathcal{X})$.
\begin{thm}
 The cone $c\text{-}\overline{{\Eff}}^1(\mathcal{X})$ is the closure of $ c\text{-}\text{\rm Big}^1(\mathcal{X})$.
\begin{proof}
First note that $c\text{-}\text{Big}(\mathcal{X}) \subseteq c\text{-}\overline{\Eff}^1(\mathcal{X})$, and  $c\text{-}\overline{\Eff}^1(\mathcal{X})$ being closed, we have $$\overline{c\text{-}\text{Big}(\mathcal{X})} \subseteq c\text{-}\overline{\Eff}^1(\mathcal{X}).$$
Now let $\alpha \in c\text{-}\overline{\Eff}^1(\mathcal{X}).$
This implies $\alpha_{X'}\geq 0$ for all $X'\in \mathcal{M}_{X}.$ Hence we conclude that $\alpha_{X'} \in \overline{c\text{-}\text{Big}(\mathcal{X})}.$  Therefore 
\begin{align*}
\alpha = \phi_{X'}(\alpha_{X'})
&\in \phi_{X'}\bigl(\overline{c\text{-}\text{Big}(X')}\bigr) \subseteq \overline{\phi_{X'}(\text{Big}(X'))}= \overline{c\text{-}\text{Big}(\mathcal{X})}
\end{align*}
This completes the proof.
\end{proof}
\end{thm}

\section{Map between $b\text{-}$divisor classes}\label{sec7}
In this section, we revisit the functoriality property of $b\text{-}$divisor classes as an application of the preceding sections. The functionalities of the $b\text{-}$ divisor classes are also discussed in \cite{BFF12}.
Let $X$ and $Y$ be two smooth projective varieties and $\mathcal{M}_X$ and $\mathcal{M}_Y$ be respective smooth models. Let $f:X\rightarrow Y$ be any smooth morphism. Now, any $Y'\in\mathcal{M}_Y$, corresponds some $X'\in\mathcal{M}_X$ such that 
\begin{equation}\label{diagram2}
\begin{tikzcd}
X' \arrow{r}{f'} \arrow[swap]{d}{\pi'_X} & Y' \arrow{d}{\pi'_Y} \\
X \arrow{r}{f} & Y
\end{tikzcd}
\end{equation}
commutes. We can consider $X' $ to be the fiber product $X\times_{Y} Y'$. Therefore, we have the pullback map $f'^*:N^1(Y')\rightarrow N^1(X')$ which induces a linear map between the $b\text{-}$ classes $\tilde{f}:c\text{-} N^1(\mathcal{Y})\rightarrow c\text{-} N^1(\mathcal{X})$, where $c\text{-}N^1(\mathcal{Y})$ and $ c\text{-} N^1(\mathcal{X})$ are Cartier $b\text{-}$ classes of co-dimension one of $Y$ and $X$, respectively. Also, observe that any commutative diagram like above corresponds to the following commutative diagram:
\begin{equation}\label{diagram3}
\begin{tikzcd}
c\text{-} N^1(\mathcal{Y}) \arrow{r}{\tilde{f}} \arrow[swap]{d} & c\text{-} N^1(\mathcal{X}) \arrow[swap]{d}\\
w\text{-} N^1(\mathcal{Y})  \arrow[swap]{d} & w\text{-} N^1(\mathcal{X}) \arrow[swap]{d}\\
N^1(Y') \arrow{r}{f'^*} & N^1(X')
\end{tikzcd}
\end{equation}
where the first vertical arrows are inclusions and the second vertical arrows are projections in the diagram. Now it is clear from the commutativity of the diagram that $\tilde{f}(\alpha)_{X'}=f'^*(\alpha_{Y'})$. In the next theorem we study the pullback of $b\text{-}$ divisor classes.

\begin{thm}
    Let $f:X\rightarrow Y$ be a smooth morphism between smooth projective varieties and $\tilde{f}:c\text{-} N^1(\mathcal{Y})\rightarrow c\text{-} N^1(\mathcal{X})$ be the corresponding linear map between the the set of Cartier $b\text{-}$ divisor classes.
    \begin{enumerate}
        \item Then $\tilde{f}\bigl(c\text{-}\Nef^1(\mathcal{Y})\bigr)\subseteq c\text{-}\Nef^1(\mathcal{X})$.
        \item If $f$ is surjective and $\tilde{f}(\alpha)$ is nef Cartier b-class, then $\alpha$ is also a nef Cartier b-class.
        
    \end{enumerate}
\end{thm}
\begin{proof}
\begin{enumerate}
    \item  If $\alpha\in c\text{-}\Nef^1(\mathcal{Y})$, then $\alpha_{Y'}$ is a nef divisor for some determination of $\alpha$, $Y'\in\mathcal{M}_Y$. Now consider the commutative diagram (\ref{diagram2}) such that $f^*(\alpha_{Y'})$ is a nef divisor on $X'$ for some $X'\in\mathcal{M}_X$. Hence $\tilde{f}(\alpha)$ is a nef Cartier $b$-class.

    \item If $\tilde{f}(\alpha)$ is a nef b-class, then there exists $X_1\in\mathcal{M}_X$ such that $\tilde{f}(\alpha)_{X_1}$ is a nef line bundle. Let $\alpha$ be determined on $Y'\in \mathcal{M}_Y$. Considering the commutative diagram (\ref{diagram2}), $\tilde{f}(\alpha)$ is also determined on $X'$. Now there exists $X''\in\mathcal{M}_X$, such that $X''\geq X',X_1$. As $f$ is a surjective map, $X''\xrightarrow{\pi'} X'\xrightarrow{f'} Y'$ where $f'$ is a surjective map followed by the birational map $\pi'$. We have $\pi'^*(f'^*(\alpha_{Y'}))=\tilde{f}(\alpha)_{X''}$ is a nef divisor. So $\alpha_{Y'}$ is a nef divisor on $Y'$ (\cite{L1}, Example 1.4.4). Hence $\alpha$ is a nef divisor.

\end{enumerate}
   
\end{proof}

\end{document}